\documentstyle[12pt]{article}
\newtheorem{defin}{Definition}
\newtheorem{prop}{Proposition}
\newtheorem{nt}{Remark}
\newtheorem{Th}{Theorem}
\newtheorem{lemma}{Lemma}
\newtheorem{cons}{Corollary}

\newfont{\sdbl}{msbm9}
\newfont{\dbl}{msbm10 at 12pt}

\newcommand{\proof}{{\bf Proof.\ }}

\newcommand{\oo}{{\cal O}}
\newcommand{\hoo}{\hat{\oo}}

\newcommand{\g}{{\cal G}}
\newcommand{\h}{{\cal H}}
\newcommand{\ad}{{\cal A}}

\newcommand{\Dim}{\mathop {\rm Dim}}
\newcommand{\Hom}{\mathop {\rm Hom}}

\newcommand{\Lim}{\mathop {\rm lim}}

\newcommand{\Spec}{\mathop {\rm Spec}}

\newcommand{\Frac}{\mathop {\rm Frac}}
\newcommand{\dm}{\mathop {\rm dim}}

\newcommand{\da}{{\mbox{\dbl A}}}
\newcommand{\dz}{{\mbox{\dbl Z}}}

\newcommand{\Z}{\dz}
\newcommand{\sdz}{{\mbox{\sdbl Z}}}
\newcommand{\sZ}{\sdz}

\newcommand{\Ker}{{\rm Ker}\:}

\newcommand{\f}{{\cal F}}
\newcommand{\lto}{\longrightarrow}

\begin{document}
\author{D.V. Osipov}

\title{Central extensions and reciprocity laws on algebraic surfaces.}
\date{}
\maketitle

\section{Introduction.}
Let
$X$
be an algebraic surface over a perfect field  $k$.
Let $C$
be an irreducible curve on $X$, $x$ a point on $C$.
From such data it is possible to construct an artinian ring $K_{x,C}$,
which is the finite direct sum of two-dimesnional local fileds
 $L_i$
(see~\cite{P}, \cite{B}, \cite{FP}).

Each two-dimensional local field
 $L_i = k_i((t_i))((u_i))$,
where the field $k_i$  is a finite extension of the field $k$.
For the point $x$, which is a smooth point on $X$,
each of fields $k_i$ coinsides with the residue field of the point $x$.
For the point $x$, which a smooth point on $X$ and $C$,
the ring $K_{x,C}$
is only a two-dimensional local field.

On a two-dimensional local field $L = k((t))((u))$
there exists a discrete valuation of rank $2$:
$$
 (\nu_1, \nu_2) \quad :  \quad  L^*  \lto  \Z \oplus \Z
$$

We define a symbol $\nu_L$:
$$
K_2 (L) \lto \Z    \qquad  \mbox{,} \qquad
\nu_{L}(f,g) =
\left|
\begin{array}{cc}
\nu_1(f) & \nu_1(g) \\
\nu_2(f) &  \nu_2(g)
\end{array}
\right|
$$
where $f, g \in L^*$.

Such symbol was used in~\cite{P2}
for the description of the intersection theory of divisors
on algebraic surface.

Also in \cite{P1} the symbol $\nu_L$
was used for the description of non-ramified
extensions of the field
$L$
for the finite field $k$.

If  $L' / L$
is non-ramified extension,
then the reciprocity map
maps
$$
 K_2(L) \ni  (f, g)  \;  \mapsto  \;
 Fr^{\nu_L(f,g)}  \in \mathop{Gal} (L' /L)
\mbox{,}
$$
where $Fr$
is the Frobenius automorphism
lifted to an element from the Galois group
 $Gal (L' / L)$.

For $K_{x,C} = \mathop{\oplus}\limits_i L_i$ we define
$$
\nu_{x,C}= \bigoplus_i \: [k_i : k] \cdot \nu_{L_i}
$$

For $\nu_{x,C}$ the following reciprocity laws hold.

{\em We fix a point $x \in X$,
then for any $f,g \in k(X)^*$
in the following sum only the finite number of terms are not equal
to zero, and
$$
\sum_{ C \ni x}   \nu_{{x,C}} (f,g) = 0     \mbox{,}
$$
where we take the sum over the all irreducible curves
 $C$ on $X$
that contain the point $x$.}

{\em
We fix an irreducible projective curve
 $C$ on a smooth surface $X$. Then
for any $f,g \in k(X)^*$ in the following sum
only the finite number of terms is not equal to zero, and
$$
\sum_{ x \in C}   \nu_{{x,C}} (f,g) = 0  \mbox{,}
$$
where we take the sum over the all points $x$
of the curve $C$.}

We interpret  $\nu_{{x,C}}$
as the commutator of lifting of elements in
some central extension of the group $k(X)^*$.
For this goal we use the notion of dimension theory,
which was
introduced by Kapranov in~\cite{Ka}.

We give the proof of the reciprocity laws
using  the splitting of central extension over
the rings of adeles connected with
the points on surfaces and with the projective curves on surfaces.

The  proof is similar to the Tate method on sum of residue
of differentials on the projective curve, see~\cite{T}, \cite{Arbar}.

During the work on this article the works~\cite{Br} and~\cite{R}
were useful for me.

I am  very grateful to A.N.~Parshin for his attention to this work
and helpful  discussions. I am grateful to H.~Kurke and A.~Zheglov
for some discussions.

\section{Construction of the group.}

\subsection{Linearly locally compact vector spaces.}

\begin{defin}{\cite[ch.III, \S 2, ex.15-21]{Bo}]}
A topological
 $k$-vector
space $V$ (over a discrete field $k$)
is called a linearly compact space,
if the following conditions hold:
\begin{enumerate}
\item the space $V$
is a complete and Hausdorff space,
\item the space $V$
has a base of neighbourhoods of
 $0$
consisting of vector subspaces,
\item
every open subspace in
 $V$
has finite codimension.
\end{enumerate}
\end{defin}

\begin{defin}{\cite[ch.2 ,\S 6]{L}}
A topological
 $k$-vector space $W$ (over a discrete field $k$)
is called a linearly locally compact space,
if it has the base of neighbourhoods of
  $0$
consisting of linearly compact open subspaces.
\end{defin}

Let a field $k'$ be a finite extension of the field $k$.
An example of linearly locally compact space is
the field of Laurent series
 $K = k'((t)))$
with the base of topology given
by subspaces
$$U_n = t^n \oo_K \mbox{,} $$
where the ring $\oo_K = k'[[t]]$
is a linearly  compact space.

\begin{lemma}
\begin{enumerate}
\item \label{a}
If a $k$-vector space $V$ is a linearly compact space,
then for any closed subspace $H$ we have that the spaces
$H$ and $V / H$ are linearly compact spaces.
\item  \label{b}
If a $k$-vector space $W$ is a
linearly locally compact space,
then for any closed subspace $E$ we have that the spaces
$E$ and $W / E$ are linearly locally compact spaces.
\end{enumerate}
\end{lemma}
\proof
Item~\ref{a} proved in \cite[ch.~2 ,\S 6]{L}.
The proof of item~\ref{b}
follows at once from item~\ref{a}
and the definition of  linearly locally compact space.
\begin{flushright}
\fbox{}
\end{flushright}

\begin{defin}
If $k$-spaces $V_1$ and $V_2$
are linearly locally compact spaces,
then ${\cal H}om (V_1, V_2)$
consists of $k$-linear  continuous
maps between
 $V_1$ and $V_2$.
\end{defin}

 The field of Laurent series
 $K=k'((t))$
is a one-dimensional local field over the field $k$.
We denote the ring of $k$-linear
continuous maps
$$ End(K/k)= {\cal H}om (K, K) \mbox{.} $$
We define the group of invertible elements of this ring by
 $$GL(K/k) = End(K/k)^*  \mbox{.} $$

\subsection{Algebra of endomorphisms of two-dimensional local field.}
Let a field $k'$ be a  finite extension of the field $k$.
We consider now the field of iterated Laurent series
$$L= k'((t))((u)) \mbox{.} $$
Such a field is called a two-dimensional local field over the field $k$.
The field $L$
has the first residue field $\bar{L} = k' ((t))$,
which is a one-dimensional local field over the field $k$.

We consider now the ring $\oo_L = k' ((t))[[u]]$.
We define for any integer $n$
a space
$$\oo_n = u^n  k'((t))[[u]] \mbox{.}$$
If $m < n$, then the  $k$-vector
space $\oo_m / \oo_n$
is a linearly locally compact space with the topology
induced by open subspaces
$ E_l = t^l u^m k'[[t,u]]  / t^l u^n k'[[t,u]]  \subset  \oo_m / \oo_n
$.

\begin{defin} \label{d1}
Let
 $L= k'((t))((u))$
be a two-dimensional local field over the field $k$.
We say that a  $k$-linear map $A : L \to L$
is from $End(L/k)$, if the following conditions hold
\begin{enumerate}
\item \label{i1}
for any integer $n$
there exists an integer $m$
such that $A \oo_{n} \subset \oo_{m} $,
\item \label{i2}
for any integer
 $m$
there exists an integer $n$
such that $A \oo_n \subset \oo_m $,
\item \label{i3}
for any integer $n_1 < n_2 $
and $m_1 < m_2$ such that $A \oo_{n_1} \subset \oo_{m_1}$
and $A \oo_{n_2} \subset \oo_{m_2}$ we have that
the induced $k$-linear map
$$   \bar{A} \quad : \quad \oo_{n_1} / \oo_{n_2} \lto
\oo_{m_1} / \oo_{m_2}
$$
belongs to
${\cal H}om
(\oo_{n_1} / \oo_{n_2} , \oo_{m_1} / \oo_{m_2}) \mbox{.} $
\end{enumerate}
\end{defin}

\begin{prop}   \label{pre}
\begin{enumerate}
\item \label{it1} $End(L/k)$ is a $k$-algebra.
\item \label{it2}
For any integer $m$
there is an embedding of algebras
$End(L/k)^{\oplus m} \hookrightarrow End(L/k)$.
\item \label{it3}
For any integer $m$
there is an embedding of algebras
$gl(m, L) \hookrightarrow End(L/k)$.
\end{enumerate}
\end{prop}
\proof
Item~\ref{it1}
follows easy from the definition of $End(L/k)$.
Therefore we give the proof of items~\ref{it2} and~\ref{it3}.
We have $L = \bar{L}((u))$.
Let $e_i$, $1 \le i \le m$
be a standart basis in $L^{\oplus m}$,
i.e.,
$$
L^{\oplus m} = \bigoplus_i  L e_i \mbox{.}
$$
We consider a $\bar{L}$-linear isomorphism
$\phi  :  L^{\oplus m}  \to L$
$$
  \phi(\sum_{i,j} a_{i,j} u^j e_i) = \sum_{i,j} a_{i,j} u^{mj+i-1},  \quad
a_{i,j} \in \bar{L} \mbox{.}
$$
Then any $k$-linear operator $A$, which acts on $L^{\oplus m}$,
under the map
\begin{equation} \label{equ1}
A \mapsto  \phi A \phi^{-1}
\end{equation}
goes to a $k$-linear operator, which acts on $L$.
In the definition~\ref{d1}
we can change the spaces $\oo_l$ (for all integer  $l$)
to the spaces $\oo_{lm}$ for any fixed $m$.
And
$$\phi (\bigoplus_i \oo_l e_i) =  \oo_{lm} \mbox{.} $$
Therefore after the easy verification we get that
the map~(\ref{equ1})
gives an embedding of algebras
$End(L/k)^{\oplus m}$
and $gl(m, L)$,
which act on $L^{\oplus m}$,
to the algebra $End(L/k)$.
\begin{flushright}
\fbox{}
\end{flushright}

From the proposition above we obtain at once the following corollary
\begin{cons}
For any natural $m$
we have an embedding of toroidal Lie algebra
$$
gl(m, k[t,u,t^{-1}, u^{-1}]) \hookrightarrow  End(L/k) \mbox{.}
$$
\end{cons}

\begin{nt}
{\em
 The field $L = k'((t))((u))= k'((u))\{\{t\}\}$,
where $k'((u))\{\{t\}\}$ denotes the following expressions:
$$
g = \sum_{i = - \infty}^{i = + \infty} a_i t^i  \quad \mbox{,} \quad
a_i \in k'((u))  \mbox{,}
$$
where
$\Lim\limits_{i \to -\infty } \nu(a_i) = 0$
and for all $i$ we have $\nu (a_i) > c_g$
for some integer $c_g$,
and  $\nu$
is the discrete valuation of the field $k'((u))$.

Then for any natural $m$  we have
a well-defined
$k'((u))$-linear isomorphism
$$ \psi  :  L^{\oplus m}  \to L
\qquad \mbox{,}  \qquad
  \psi(\sum_{i,j} a_{i,j} t^j e_i) = \sum_{i,j} a_{i,j} t^{mj+i-1},  \quad
a_{i,j} \in k((u)) \mbox{.}
$$

 The isomorphism $\psi$
gives an another embedding of algebra $gl(m, L)$ to the algebra $End(L/k)$
by means of the action of $gl(m, L)$ on $L^{\oplus m}$.
}
\end{nt}

Further we will need the following corollary
\begin{cons} \label{sl}
An action  of the field $L$ on $L$ by multiplications
gives a canonical embedding
$$
L  \hookrightarrow End(L/k)  \mbox{.}
$$
\end{cons}

\begin{nt} {\em
The algebra $End(L/k)$
was defined by A.~Beilinson in~\cite{B}
for a mu\-lti\-di\-men\-si\-o\-nal local field.
The field
$$
L =
\mathop{\Lim_{\rightarrow}}_n \mathop{\Lim_{\leftarrow}}_{m \ge n}
\oo_n / \oo_m   \mbox{.}
$$
}
\end{nt}
It is not difficult to understand
that the definition of $End(L/k)$
copies the definition of morphisms in $Ind Pro$-category.
See the abstract description of this in appendix to~\cite{B1}.

\subsection{Topology of two-dimensional local field.}

On the field $L= k'((t))((u))$
there is a standard topology,
where the base of neighbourhoods of zero are
 $k'$-subspaces
\begin{equation} \label{top}
 \sum_i U_{l_i} u^i  + \oo_m \mbox{,}
\end{equation}
where $U_{l_i} = t^{l_i} k'[[t]]$
are open  $k'$-vector
subspaces from $k'((t))$ and $m$
is an integer number.
We note, that the spaces
 $\oo_m$
are closed in this topology,
and the  subspaces $\oo_m / \oo_n$
are linearly locally compact spaces with the induced topology.

\begin{lemma}  \label{och}
The algebra $End (L/k)$
acts continuously
on the field $L$.
\end{lemma}
\proof
We consider any element $A \in End (L/k)$
and some open subspace
$V = \sum_i U_{l_i} u^i  + \oo_m $.
Then from the definition of the algebra $End (L/k)$
it follows
that there exists an integer $n$
such that $A \oo_n \subset \oo_m$.
From the definition $End (L/k)$ it follows that
there exists $m_1 < m$ such that
$A \oo_{n-1} \subset \oo_{m_1}$.
There exists an open subspace
$E_{k_1} \subset \oo_{n-1} / \oo_{n}$
such that $A E_{k_1} \subset (V \cap \oo_{m_1}  / V \cap \oo_m)$.
Now there exists $m_2 < m_1$ such that
 $A \oo_{n-2} \subset \oo_{m_2}$.
Then there exists an open subspace
 $E_{k_2} \subset (\oo_{n-2} / \oo_{n}) $
such that
$$A E_{k_2} \subset (V \cap \oo_{m_2}  / V \cap \oo_m) \qquad
\mbox{and} \qquad E_{k_2} \cap (\oo_{n-1} / \oo_{n})   \subset E_{k_1}
\mbox{.}$$
In the same way we define now by induction for any natural $i$
an integer $m_i < m_{i-1}$ such that $A \oo_{n-i} \subset \oo_{m_i}$,
and we define an open subspace
$E_{k_i} \subset \oo_{n-i} / \oo_n$
such that
$$A E_{k_i} \subset (V \cap \oo_{m_i}  / V \cap \oo_m ) \qquad \mbox{and}
\qquad  E_{k_i} \cap (\oo_{n-i+1} / \oo_{n})   \subset E_{k_{i-1}} \mbox{.}$$
Then the space $W = \sum_i E_{k_i} + \oo_n$
is open in $L$, and we have  $A W \subset V$.
\begin{flushright}
\fbox{}
\end{flushright}

\begin{nt} {\em
Let  $X/k$
be an algebraic surface over a field $k$.
For any irreducible curve $C \subset X$
and any point $x \in C$
one can canonically define a ring
$$
K_{x,C} = k(X) \cdot
\mathop{\Lim_{\leftarrow}}_m
({\oo_{x, X}}_{(J_C)}  \cdot \hat{\oo}_{x,X}) / J_C^m
\mbox{,}
$$
where $J_C$
is the ideal  of the curve $C$
in the local ring $\oo_{x, X}$
of the point $x$ on the surface $X$,
the ring ${\oo_{x, X}}_{(J_C)}$
is the localization
 of the ring $\oo_{x,X}$
along the prime ideal $J_C$.

Then the ring
$$
K_{x,C} = \oplus K_i  \mbox{,}
$$
where every $K_i$ is a two-dimensional local field over the field $k$,
and the direct sum is  over the finite set.
On the ring $K_{x,C}$
there is a natural topology of inductive and projective limits,
which comes from the construction of the ring.
On the field $K_i$
we have the induced topology.
In each two-dimensional local field $K_i$
we can choose
the system of local parameters
 $t_i, u_i$ and isomorphism
$k_i((t_i))((u_i)) = K_i$
such that the topology on the field $K_i$
coincides with the standart topology
(\ref{top})
of the field $k_i((t_i))((u_i))$ (see~\cite{Y}).
(For the smooth point $x$ on $X$ each of fields $k_i$
coinsides with the residue field of the point $x$.)
Another choice of local parameters gives a continuous automorphism
of the field $K_i$,
which preserves the subrings $\oo_m$ in  $K_i$.
Hence we have that the algebras $End(K_i/k)$
do not depend on the choice of local parameters $t_i, u_i$.
}
\end{nt}

\begin{defin}
Let  $L = k'((t))((u))$
be a two-dimensional local field over the field $k$.
Then we define a group
$$
GL(L/k) = End(L/k)^*  \mbox{.}
$$
\end{defin}

\section{Central extension.}

\subsection{Dimension theory.}

Let
 $V$
be a linearly locally compact space over the field $k$.
For any two open linearly compact  $k$-subspaces
$A$ and $B$ from $V$ we define
$$
 [A \mid B]_1   \quad =  \quad \dm\nolimits_k \frac{B}{A \cap B} -
\dm\nolimits_k \frac{A}{A \cap B}
$$

This definition is correct, since every open subspace
in a linearly compact vector space has finite codimension.
We have the following property.

\begin{lemma}
For any three open linearly compact $k$-subspaces of a space
$V$ the following formula holds
$$
[A \mid B]_1  + [B \mid C]_1 =  [A \mid C]_1  \mbox{.}
$$
\end{lemma}
\proof
It is easy to understand that for any open linearly compact
$k$-subspace $U \subset A$, $U \subset B$ the following formula holds
$$
 [A \mid B]_1   \quad =  \quad \dm\nolimits_k B/U -
\dm\nolimits_k A/U   \mbox{.}
$$
Now we choose such an $U \subset A $, $U \subset B$, $U \subset C$,
and we obtain the statement of lemma.
\begin{flushright}
\fbox{}
\end{flushright}

\begin{defin}{\cite{Ka}}
Let
 $V$
be a linearly locally compact space over the field $k$.
An integer value  function $d$ on the set
of open linearly compact subspaces of the space $V$
is called the dimension theory,
if for any two open linearly compact subspaces
 $A$ and $B$ from $V$
we have
$$
d(B) = d(A) + [A \mid B]_1 \mbox{.}
$$
\end{defin}

By $\Dim (V)$  we will denote the set of dimension
theories of the space $V$.

The following proposition is from~\cite{Ka}.
\begin{prop} \label{Kap}
\begin{enumerate}
\item $\Dim (V)$ is a $\Z$-torsor.
\item
For each short exact sequence of $k$-spaces
\begin{equation} \label{tri}
0 \lto V_1 \lto V \lto  V_2  \lto 0  \mbox{,}
\end{equation}
where
 $V$
is linearly locally  compact,
  $V_1$ is closed in $V$,
and  $V_1$
is with the restriction topology,
 $V_2$ is with the factor topology,
we have the natural  isomorphism of   $\Z$-torsors
$$
\Dim (V')  \otimes_{\sZ} \Dim(V'') \lto \Dim(V)  {.}
$$
\end{enumerate}
\end{prop}
\proof
Let
 $d$
be a dimension theory on $V$,
then after the addition of an integer number to the values of $d$ we obtain
the new dimension theory on $V$.
This action of $\Z$  makes
from  $\Dim (V)$ a $\Z$-torsor.

Let $d_1$
be a dimension theory on $V_1$,
and $d_2$
be a dimension theory on $V_2$,
then for any linearly compact subspace
$U$ from $V$
we define
$$
d_1 \otimes d_2 (U) = d_1(U \cap V_1) + d_2(U / U \cap V_1) \mbox{.}
$$
Now $d_1 \otimes d_2$    is a well-defined dimension theory
on  $V$.
\begin{flushright}
\fbox{}
\end{flushright}

\begin{nt}
{\em
Any linearly compact subspace
 $W$ from $V$
determines canonically
 $d_W  \in \Dim(V)$ by the following rule
$$
d_{W} (U) = [W \mid U]_1 \mbox{.}
$$
 }
\end{nt}

\subsection{$\Z$-torsor
of relative dimension theories.}  \label{razdel}
We consider a $k$-vector space
$$M = \prod_{i \in I} L_i \mbox{,}  $$
where $I$
is some set, and every
$L_i= k_i((t_i))((u_i))$
is a two-dimensional local field over the field $k$.
On each $L_i$ there is the topology,
and therefore on $M$ there is the topology of product.

For any two closed $k$-subspaces $W_1$,  $W_2$ from $M$
such that there exists a closed $k$-subspace $W \subset W_1$,
$W \subset W_2$ and
the spaces  $W_1 /W$ and $W_2 / W$ are
linearly locally compact with induced topology from $M$
we define a $\Z$-torsor
$$
[W_1 \mid W_2 ; W  ]_2 = \Hom\nolimits_{\sZ} (\Dim(W_1 / W) , \Dim(W_2 / W)).
$$

We have the following evident lemma:
\begin{lemma} \label{oche}
Let $W_1, W_2, W_3, W$ are closed $k$-subspaces from $M$
such that $W_1 / W$, $W_2 /W$, $W_3 /W$
are linearly locally compact. Then there is a canonical isomorphism
of  $\Z$-torsors
$$
[W_1 \mid W_2 ; W  ]_2  \otimes_{\sZ} [W_2 \mid  W_3 ; W]_2 \lto
[W_1 \mid  W_3 ; W]_2
$$
\end{lemma}

\begin{lemma}
Let a $k$-subspace $W' \subset W$
satisfies the same conditions as $W$.
Then there exists a canonical isomorphism
$$
[W_1 \mid  W_2 ; W']_2  \lto [W_1 \mid  W_2 ; W]_2
$$
\end{lemma}
\proof
From proposition~\ref{Kap}
we have
$$
\Dim(W_1 / W')= \Dim (W_1 / W) \otimes_{\sZ} \Dim(W / W')
$$
$$
\Dim(W_2 / W')= \Dim (W_2 / W) \otimes_{\sZ} \Dim(W/W') \mbox{.}
$$
Therefore let
$\phi \in \Hom\nolimits_{\sZ} (\Dim(W_1 / W) , \Dim(W_2 / W))$,
then for a fixed $a \in \Dim (W / W')$
we define $\phi' (d \otimes a) = \phi (d) \otimes a$,
where $d \in  \Dim (W_1 /W)$.
Then
$$
\phi' \in \Hom\nolimits_{\sZ} (\Dim(W_1 / W') , \Dim(W_2 / W'))  \mbox{,}
$$
and an isomorphism $\phi \mapsto \phi'$ does not depend on the choice of $a$.
\begin{flushright}
\fbox{}
\end{flushright}

Due to the last lemma we can give the following definition.
\begin{defin}
Let $k$-subspaces $W_1$ and $W_2$ be as above. We define
$$
[W_1 \mid W_2]_2 \quad  =
\quad \mathop{\Lim_{\longleftarrow}}_W \quad  [W_1 \mid W_2 ; W]_2
$$
\end{defin}

\begin{prop}  \label{predlozh3}
\begin{enumerate}
\item $[W_1 \mid W_2]_2$ is a $\Z$-torsor.
\item
Suppose that for closed
 $k$-subspaces $W_1, W_2, W_3$
there exists a closed $k$-subspace
$W \subset W_1$, $W \subset W_2$, $W \subset W_3$
such that $W_1 / W$, $W_2 / W$ ¨ $W_3 / W$
are linearly locally compact spaces,
then there exists a canonical isomorphism
$$
[W_1 \mid W_2]_2 \otimes_{\sZ} [W_2 \mid W_3]_2 \lto [W_1 \mid W_3]_2 \mbox{,}
$$
and for any four closed $k$-subspaces $W_1, W_2, W_3, W_4$
such that there exists a closed
 $k$-subspace $W$
with linearly locally compact spaces
  $W_1 / W$, $W_2 / W$, $W_3 / W$, $W_4 / W$
we have that the following diagram is commutative.
\begin{equation} \label{ass1}
\begin{array}{ccc}
[W_1 \mid W_2]_2
\otimes_{\sZ} [W_2 \mid W_3]_2
\otimes_{\sZ} [W_3 \mid W_4]_2
&
\lto
&
[W_1 \mid W_2]_2 \otimes_{\sZ} [W_2 \mid W_4]_2
\\
\downarrow & & \downarrow
\\
{[W_1 \mid W_3]_2
 \otimes_{\sZ}
[W_3 \mid W_4]_2}
&
\lto
&
[W_1 \mid W_4]_2  \mbox{.}   \\
\end{array}
\end{equation}
\end{enumerate}
\end{prop}
\proof
The sets $[W_1 \mid W_2 ; W]_2$  are $\Z$-torsors.
We  identify all these sets by means of projective limit
with isomorphic maps.
Therefore $[W_1 \mid W_2]_2$
is a $\Z$-torsor as well.

We check at once by lemma~\ref{oche} the other statements for the fixed $W$,
and they commute with the replacement
  $W$ by $W'$.
Therefore we have these statements
for
 $[W_1 \mid W_2]_2$
 after the passing
to projective limit.
\begin{flushright}
\fbox{}
\end{flushright}

\subsection{Central extension $\widehat{GL}_W$.} \label{r}

Now let
 $M = L =k'((t))((u))$
be a two-dimensional local field over the field $k$.
We recall that  $\oo_n = u^n k'((t))[[u]]$.
We consider a closed $k$-subspace
$W \subset L$
such that there exist  integers $n_1 \le n_2$ with the property
\begin{equation} \label{usl}
 \oo_{n_1} \subset W \subset \oo_{n_2} \mbox{.}
\end{equation}

Then for any element
 $g \in GL(L/k)$
we get from definition~\ref{d1}
that there exist integers
 $k \le l $ such that
$$
\oo_k \subset  gW  \subset  \oo_l  \mbox{.}
$$
Therefore for any
 $l_1 \ge l$  a $k$-vector space
$gW / \oo_{l_1}$
is a
linearly locally compact space.
Therefore for any
 $l_2 \ge \max{(l_1, n_2)}$
a  $\Z$-torsor  $[W \mid gW ; \oo_{l_2}]_2$
is well-defined,
and  therefore
a $\Z$-torsor $[W \mid gW]_2$
 is well-defined as well.

Moreover, for any
 $h_1, h_2 \in GL(L/k)$
a $\Z$-torsor  $[h_1 W \mid h_2 W]_2$ is well-defined.

By lemma~\ref{och} the group $GL(L/k)$ acts   on $L$ continuously.
Therefore for any two elements
 $g_1, g_2 \in GL(L/k)$
we have a canonical isomorphism
$$g_1 [h_1 W \mid h_2 W]_2 \lto [g_1 h_1 W \mid g_1 h_2 W]_2   $$
such that
\begin{equation} \label{ass2}
g_2 (g_1 (d)) = (g_2 g_1) (d) \mbox{,}
\end{equation}
for any $d \in [h_1 W \mid h_2 W]_2$.

\begin{defin}   \label{def8}
We define a set of pairs
$$
\widehat{GL}_W = (g, d)  \mbox{,}
$$
where $g \in GL(L/k)$, $d \in [W \mid gW]_2$.
\end{defin}

We can multiply these pairs by a rule
\begin{equation} \label{um}
(g_1, d_1) (g_2, d_2) = (g_1 g_2, d_1 g_1 (d_2))  \mbox{,}
\end{equation}
here we used an isomorphism
$ [W \mid g_1 W]_2 \otimes_{\sZ} [g_1 W \mid g_1 g_2 W]_2
\to [W \mid g_1 g_2 W]_2 $,
and also that
 $g_1 (d_2) \in [g_1 W \mid g_1 g_2 W]_2$.

\begin{prop} \label{predlo}
\begin{enumerate}
\item
The multiplication, which was defined by the equality~(\ref{um}),
makes the group from the set $\widehat{GL}_{W}$.
\item  \label{pu}
For any other closed $k$-subspace
$W'$ such that
$\oo_{m_1} \subset  W'  \subset \oo_{m_2}$
for some integers
 $m_1$, $m_2$
we have a canonical isomorphism
$$
\widehat{GL}_{W'} \lto \widehat{GL}_{W} \mbox{.}
$$
\item
We have a central extension
\begin{equation} \label{centra}
0 \lto \Z \lto \widehat{GL}_W \stackrel{\pi}{\lto} GL(L/k) \lto 1 \mbox{.}
\end{equation}
\item  \label{pun}  Let $H$ be a subgroup in $GL(L/k)$
such that for any element $h \in H$ we have $hW \subset W$,
then the central extension~(\ref{centra})
splits over the subgroup $H$.
\end{enumerate}
\end{prop}
\proof
It is clear that the unit of the group
is an element $(e, d_0)$,
where $d_0$
is a function on the zero space with zero value.
Also from the construction it follows the existence
of the inverse element.
The associativity of multiplication
follows from diagram~(\ref{ass1}) and equality~(\ref{ass2}).
 Therefore $\widehat{GL}_W$ is a group.

Now we consider the subspace $W'$.
We fix  an arbitrary element  $ c \in [W' \mid W]_2$.
Then a map
$$
(g,d)  \mapsto (g, c d g(c^{-1}))   \quad \mbox{,} \qquad
 d \in [W \mid gW]_2  \mbox{,}
$$
gives an isomorphism of  $\widehat{GL}_W$ onto $\widehat{GL}_{W'}$.
It is easy to see that this isomorphism does not depend
on the choice of the element $c \in [W' \mid W]_2$.

A map
$$
(g,d)  \mapsto g
$$
gives the map
 $\pi$ of the group $\widehat{GL}_W$ onto the group $GL(L/k)$
with the kernel $\Z$.

Item~\ref{pun}
follows from the definition,
since $[hW \mid W]_2 = [W \mid W]_2$
is a trivial $\Z$-torsor.
And we can take the zero element as a section.
\begin{flushright}
\fbox{}
\end{flushright}

\subsection{Symbol $\nu_L$.}
Let $L = k'((t))((u))$
be a two-dimensional local field over the field $k$,
$\bar L = k'((t))$
is its residue field.
The field $L$
has a discrete valuation of rank~$2$:
 $$(\nu_1, \nu_2) \quad : \quad  L^* \to \dz \oplus \dz  \mbox{,}$$
where  $\nu_2$
is the usual discrete valuation with respect to the local parameter
 $u$,
and
$$
\nu_1 (b) = \nu_{\bar{L}} ( \overline{b u^{-\nu_2(b)}})  \mbox{,}
$$
where $\nu_{\bar{L}}$
is the discrete valuation of the field $\bar{L}$.
We note that the valuation
$\nu_2$
is determined canonically by $L$,
but  $\nu_1$
depends on the choice of local parameter
 $u$ in the field $L$.

We define a map
$$
\nu_L \quad : \quad L^* \times L^*  \lto \Z
$$
by the following formula
\begin{equation} \label{yav}
\nu_{L}(f,g) =
\left|
\begin{array}{cc}
\nu_1(f) & \nu_1(g) \\
\nu_2(f) &  \nu_2(g)
\end{array}
\right|
\end{equation}

\begin{prop}
\begin{enumerate}
\item The map  $\nu_L$
does not depend on the choice of local parameters
$t,u$ of two-dimensional local field $L$.
\item The map $\nu_L$
is a bimultiplicative and skew-symmetric map.
\end{enumerate}
\end{prop}
\proof
The proof of this proposition follows from
explicit formula~(\ref{yav}).
\begin{flushright}
\fbox{}
\end{flushright}

\begin{nt} {\em
One can demonstrate that
the map $\nu_L$
coincides with the composition of the following maps
$$
L^* \times L^*  \lto K_2(L) \stackrel{\partial_2}{\lto} \bar{L}^*
\stackrel{\partial_1}{\lto} \dz \mbox{,}
$$
where the map $\partial_i$ is
the boundary map in algebraic $K$-theory.
The map $\partial_2$
coincides with the tame symbol with respect to the
discrete valuation $\nu_2$ of the field $L$.
The map $\partial_1$ coincides with the discrete valuation
 $\nu_{\bar{L}}$ of the field $\bar{L}$.
}
\end{nt}

\begin{nt}{\em
The symbol $\nu_L$
was used in~\cite{P2}, \cite{Lom}
for the description of intersection index of divisors on algebraic
surface by means of adelic ring of algebraic surface.

Besides, we note that the symbol
 $\nu_L$
appears
for the description
of non-ramified extensions of two-dimensional local fields
when the field
 $k$ is finite, see~\cite{P}.
}
\end{nt}

\subsection{Commutator $<\; , \;>_L$.}
We consider the central extension~(\ref{centra})
$$
0 \lto \Z \lto \widehat{GL}_W \stackrel{\pi}{\lto} GL(L/k) \lto 1 \mbox{.}
$$

\begin{defin}  \label{pod}
For any two elements $a, b \in GL(L/k)$
such that the commutator $[a,b]=1$  we define
an element
$$<a, b>_L = [a', b']  \in \Z \mbox{,}  $$
where $a' \in G_W$ and $b' \in G_W$
are any under condition $\pi (a') = a$,
$\pi (b') = b$.
\end{defin}

It is clear that
 $<a,b>_L$
does not depend on the choice of corresponding
 $a'$
and $b'$. Besides, $<a,b>_L$
does not depend on the choice of $W$,
since the central extension~(\ref{centra})
does not depend on the choice of $W$
due to item~\ref{pu}
of proposition~\ref{predlo}.

\begin{prop} \label{kocy}
\begin{enumerate}
\item The map $<\; , \; >_L$
is a  bimultiplicative and  skew-symmetric map.
\item
\begin{equation}  \label{koc}
<a , b >_L = c(a,b) - c(b,a)  \mbox{,}
\end{equation}
where $c(\;, \;)$
is any cocycle
which gives the central extension~(\ref{centra}).
\end{enumerate}
\end{prop}
\proof
The skew-symmetric property is clear.
We prove the bimultiplicativic property.
For any group we have the following formula
$$
[fg, h] = [f,[g,h]][g,h][f,h]  \mbox{.}
$$
In our case $[b', c']$
is a central element, and therefore
$$
<ab,c>_L = <a,c>_L + <b,c>_L       \mbox{.}
$$
Formula~(\ref{koc})
follows from the explicit construction of a cocycle by the central extension
and the definition of $<\;, \; >_L$.
\begin{flushright}
\fbox{}
\end{flushright}

We recall that $L = k'((t))((u))$.
From corollary~\ref{sl} of proposition~\ref{pre} we have that
 $L^*  \subset GL(L/k)$.
\begin{Th}   \label{te1}
For any elements
$f, g  \in L^* $  we have
$$
<f,g>_L= [k' : k]  \cdot \nu_L(f,g)
$$
\end{Th}
\proof
Let $W = \oo_L = k'((t))[[u]]$.
Both $<\; ,\; >_L$ and $[k' : k]  \cdot  \nu_L$
are bimultiplicative and skew-symmetric.
We have a multiplicative decomposition
$$
L^* = \bar{L}^* \times u^{\sdz} \times {\cal U}^1_L  \mbox{,}
$$
where
${\cal U}^1_L = 1 + u k'((t))[[u]]$.
Therefore for the proof of theorem it is enough to consider the following
cases.
\begin{enumerate}
\item Let $f,g \in \bar{L}^* \times {\cal U}^1_L$.
Then $\nu_L (f,g) = 0$.
On the other hand
 $f \oo_L = \oo_L $, $g \oo_L = \oo_L$.
Therefore from item~\ref{pun}  of proposition~\ref{predlo}
it follows that $ <f,g>_L = 0 $.
\item Let $f \in \bar{L}^*$, $g = u^{-1}$.
Then $\nu_L (f,u^{-1}) = - \nu_{\bar{L}} (f)$.
We fix any element $d$ from
$$  \Dim(\bar{L}) = \Dim (\oo_L / u^{-1} \oo_L) =
\Hom\nolimits_{\sdz} (\Dim (0), \Dim (\oo_L / u^{-1} \oo_L)) =
[\oo_L \mid u^{-1} \oo_L]_2 \mbox{.}$$
Let $f' = (f ,d_0)$, $g' = (g, d)$
be elements from $\widehat{GL}_W$,
where $d_0$ is the zero function.
Let $k'$-subspace $U = k'[[t]]$.
Then
$$f' g' =  (f u^{-1}, f (d) )
\qquad \mbox{and} \qquad g' f' = (f u^{-1}, d) \mbox{.}$$
Therefore
$$ <f,g>_L =  [f', g'] =
 f (d) (U) - d(U) =
 d( f^{-1}(U)) -d (U)= [U \mid f^{-1} U]_1 =
 - \nu_{\bar{L}}(f) \cdot [k' : k] \mbox{.}$$
\end{enumerate}
\begin{flushright}
\fbox{}
\end{flushright}

\subsection{Additivity
$<\; , \;>_M = <\; , \;>_{M_1} + <\; , \;>_{M_2}$.} \label{s3.6}
 Similar to section~\ref{razdel} we consider  a $k$-space
$$M = \prod_{i \in I} L_i \mbox{,} $$
where every $L_i$
is a two-dimensional local field over the field $k$.
We fix a $k$-subspace $W_i \subset L_i$ for all $i \in I$
such that for every
 $W_i$
condition~(\ref{usl}) from section~\ref{r} holds.
We consider a group $H$ such that
for any $i$
we have an embedding
$$
H  \subset GL(L_i/k)  \mbox{,}
$$
and therefore we have an action of the group $H$
on each $L_i$.
We suppose that for any element
 $h \in H$   for almost all $i \subset I$
we have
$$
h W_i \in W_i  \mbox{.}
$$

Let $$W = \prod_{i \in I}  W_i$$  be a
$k$-subspace in $M$.
Then similar to definition~\ref{def8}      we define
a group $\widehat{H}_W $
as the set of pairs $(h,d)$,
where $h \in H$, $d \in [W \mid gW]_2$.
We have the following central extension
$$
0 \lto \Z \lto \widehat{H}_W \lto H \lto 1  \mbox{.}
$$
For any two elements $a , b \in H$
such that $[a,b] =1$
we can like in definition~\ref{pod}
define $<a,b>_M = [a', b']$,
where $a'$ and $b'$
are arbitrary liftings
of elements $a$ and $b$
to
$
\widehat{H}_W
$.

\begin{prop}   \label{predl7}
Let
$I = I_1 \cup I_2$
be the set of indices.
We  consider
$k$-spaces
$$ M_l = \prod_{i \in I_l} L_i  \qquad \mbox{,} \qquad
W_l = \prod_{i \in I_l} W_i \qquad  \mbox{,}  \qquad
 l = 1,2  \mbox{.} $$
Then for any commuting  elements
 $a,b \in H$ we have
\begin{equation} \label{sum}
<a,b>_M = <a,b>_{M_1}  + <a,b>_{M_2}  \mbox{.}
\end{equation}
\end{prop}
\proof We have
$$  M = M_1 \times M_2   \qquad
  \mbox{and} \qquad W = W_1 \times W_2 \mbox{.}$$
For any $h \in H$ we have
$$
[hW \mid W]_2 = [h W_1 \mid W_1]_2 \otimes_{\sZ} [h W_2 \mid W_2]_2 \mbox{.}
$$
Therefore the central extension
$$
0 \lto \Z \lto \widehat{H}_W \lto H \lto 1  \mbox{.}
$$
is obtained as the sum of central extensions
$$
0 \lto \Z \lto \widehat{H}_{W_1} \lto H \lto 1  \qquad  \mbox{and}
$$
$$
0 \lto \Z \lto \widehat{H}_{W_2} \lto H \lto 1
$$
and the addition of kernel $\Z \oplus \Z \lto \Z$.
Therefore formula~(\ref{sum})
follows from  formula~(\ref{koc}).
\begin{flushright}
$  \fbox{}  $
\end{flushright}

\section{Reciprocity laws around points.}
Let $X$
be an algebraic surface over a field $k$.
We fix a smooth point $x$ on $X$.
Let $\hat{\oo}_{x, X}$
be the completed local  ring of the point $x$ on the surface $X$.
We consider an adelic ring of the point $x$:
$$
\da_x = \left\{ \{f_{F} \} \in \prod_{F \subset \hoo_{x,X}} K_{x,F} \quad
\mbox{such that} \quad f_F \in \oo_{K_{x,F}}
 \quad
 \mbox{for almost all}
  \quad
  F
\mbox{,}
\right\}
$$
where the product is over all prime ideals
 $F$
  of the ring  $\hat{\oo}_{x,X}$, which have codimension  $1$,
    and
$$
K_{x,F} = \Frac \; ( {\mathop {\Lim_{\longleftarrow}}_m}
 \mathop{\hoo_{x,X}}\nolimits_{ (F  ) }  / F^m
\mathop{\hoo_{x,X}}\nolimits_{ (F   ) })  \mbox{.} $$

The ring $\mathop{\hoo_{x,X}}_{(F)}$
is a localization of the ring
 $\hoo_{x,X}$
along the prime ideal $(F)$.

Let $\f \subset \Frac(\hat{\oo}_{x,X})$   be a free
 $\hat{\oo}_{x,X}$-module.
We consider a complex
${\ad_x}(\f)$:
$$
\Frac(\hat{\oo_x}) \; \times \; \prod_{F \subset \hoo_{x,X}}
(\oo_{K_{x,F}} \otimes_{\hat{\oo}_{x,X}} \f) \lto \da_x \mbox{.}
$$

\begin{prop}  \label{proposi}
The cohomology groups of the complex
 ${\ad_x}(\f)$
are linearly locally compact
 $k$-spaces.
\end{prop}
\proof
The module $\f$
gives a locally free sheaf on a scheme
 $\Spec \hat{\oo}_{x,X}$.
This sheaf we will denote by the same letter $\f$.
Then from the definition of adeles on arbitrary schemes
(see~\cite{B}, \cite{H}) we get that
the complex
 ${\ad_x}(\f)$
is an adelic complex of the sheaf $\f$
on a one-dimensional scheme
$$U = \Spec \hat{\oo_x}  \setminus x$$
($U$
is a formal neithbourhood of the point without the point.)
Therefore the cohomology groups of the complex ${\ad_x}(\f)$
coincide with the cohomology groups $H^* (U,  \f \mid_U)$.

But for any $i$  the following formula holds
\begin{equation}    \label{exty}
H^i (U, \f \mid_U) = \mathop{\mathop{\lim}_{\lto}}_n \,
{\rm Ext}^i_{\hat{\oo}_{x,X}} (m_x^n \, , \, \f) \mbox{,}
\end{equation}
where $m_x$
is the maximal ideal of the ring $\hat{\oo}_{x,X}$.

Equality~(\ref{exty}) follows from the following formula
 $$
H^0 (U , \g \mid_U) =
\mathop{\mathop{\lim}\limits_{\lto}}\limits_n
\Hom\nolimits_{\hat{\oo}_{x,X}}
(m_x^n,  \g)
\mbox{,} $$
where $\g$
is any module over the ring
 $\hat{\oo}_{x,X}$
(see~\cite[ch.3, ex.~3.7]{Ha}).

Now we apply a functor $\Hom(\cdot , \f)$ to an exact sequence  of
 $\hat{\oo}_{x,X}$-modules
\begin{equation} \label{tochn}
0  \lto m_x^n  \lto \hat{\oo}_{x,X} \lto \hat{\oo}_{x,X}/m_x^n  \lto 0
\end{equation}
We get an exact sequence
$$
0 \lto \f \lto
\Hom\nolimits_{\hat{\oo}_{x,X}} (m_x^n, \f)
\lto
{\rm Ext}^1_{\hat{\oo}_{x,X}}
(\hat{\oo}_{x,X} / m_x^n \, , \, \f)
 \lto 0
$$

We have always that $k$-spaces
${\rm Ext}^i_{\hat{\oo}_{x,X}}
(\hat{\oo}_{x,X} / m_x^n \, , \, \f) $
are finite dimensional  (see lemma~\ref{lk} below).
Therefore $H^0 (U, \f \mid_U)$
is a locally linear compact space,
where a linearly compact subspace is
 $\f$
with a base of topology given by the powers of maximal ideals
 $m_x^n \f$.

For the smooth point
 $x$ we have
$
{\rm Ext}^1_{\hat{\oo}_{x,X}}
 ( \hat{\oo}_{x,X}/m_x^n \, , \, \f)  = 0
$.
Therefore for the smooth point $x$ we have
 $H^0 (U, \f \mid_U) = \f$.

From exact sequence~(\ref{tochn})
we get that
$$
{\rm Ext}^1_{\hat{\oo}_{x,X}}
 ( m_x^n \, , \, \f)  =
{\rm Ext}^2_{\hat{\oo}_{x,X}}
 ( \hat{\oo}_{x,X}/m_x^n \, , \, \f)  \mbox{.}
$$
Therefore
$$
H^1 (U, \f \mid_U) =
\mathop{\mathop{\lim}_{\lto}}_n \,
{\rm Ext}^2_{\hat{\oo}_{x,X}} (\hat{\oo}_x / m_x^n \, , \, \f) \mbox{.}
$$
is an inductive limit of finite dimensional
 $k$-spaces, i.e.,
a discrete space, which is a linearly  locally compact space.
\begin{flushright}
$  \fbox{}  $
\end{flushright}

After proposition~\ref{proposi}
we can define a $\Z$-torsor

$$
\Dim ({\ad_x}(\f)) = \Hom\nolimits_{\sdz}
(\Dim (H^1 ({\ad_x}(\f) )) \; , \; \Dim ( H^0 ( {\ad_x}(\f)) ))
$$

Let $\f \subset \Frac(\hat{\oo}_{x,X})$
be a free
 $\hat{\oo}_{x,X}$-module.
We define a
 $k$-subspace
$$
W_{\f} =
\prod_{F \subset \hoo_{x,X}}  \oo_{K_{x,F}}
 \otimes_{\hat{\oo}_{x,X}} \f
\quad  \subset  \quad  \ad_x
$$

Since $\ad_x  \subset \prod\limits_{F \subset \hoo_{x,X}}  K_{x,F}$,
we can define the following
$\Z$-torsor
similar to section~\ref{razdel}.

\begin{defin}
For any two free
 $\oo_{x,X}$-modules $\f$ and $\g$ such that
$\f \subset \Frac(\hat{\oo}_{x,X})$
and $\g \subset \Frac(\hat{\oo}_{x,X})$
we define
$$
[W_{\f} \mid W_{\g}]_2 =
\quad \mathop{\Lim_{\longleftarrow}}_{W_{\h}} \quad  [W_{\f} \mid W_{\g} \;
; \; W_{\h}]_2
\mbox{,}
$$
where the limit is taken over all free  $\oo_{x,X}$-modules
$\h$ such that $\h \subset \f$ and $\h \subset \g$.
\end{defin}

We recall that from proposition~\ref{predlozh3}
we have a canonical isomorphism of $\Z$-torsors:
$$
[W_{\f} \mid W_{\g}]_2 \otimes_{\sZ} [W_{\g} \mid W_{\h} ]_2
\lto [W_{\f}  \mid W_{h}]_2 \mbox{,}
$$
where $\f$, $\g$ and $\h$
are any three free
 $\hoo_{x,X}$-modules
such that
$$\f \subset \Frac(\hat{\oo}_{x,X})   \qquad  \mbox{,} \qquad
\g \subset \Frac(\hat{\oo}_{x,X}) \qquad  \mbox{,} \qquad
\h \subset \Frac(\hat{\oo}_{x,X})  \mbox{.}
$$

We have the following proposition
\begin{prop}  \label{stand}
For any two free
 $\oo_{x,X}$-modules $\f$ and  $\g$
such that
$$\f \subset \Frac(\hat{\oo}_{x,X}) \qquad \mbox{and}
\qquad   \g \subset \Frac(\hat{\oo}_{x,X})  $$
we have the following canonical isomorphism of
$\Z$-torsors
\begin{equation} \label{ttt}
[W_{\f} \mid W_{\g}]_2
\lto
\Hom\nolimits_{\sZ} ( \Dim ({\ad_x}(\f)), \Dim ({\ad_x}(\g)) )
\end{equation}
and for any three free
$\oo_{x,X}$-modules $\f$, $\g$ and $\h$
such that
$$\f \subset \Frac(\hat{\oo}_{x,X})   \qquad  \mbox{,} \qquad
\g \subset \Frac(\hat{\oo}_{x,X}) \qquad  \mbox{,} \qquad
\h \subset \Frac(\hat{\oo}_{x,X})
$$
we have that isomorphism~(\ref{ttt})
maps an isomorphism
$$
[W_{\f} \mid W_{\g}]_2
\otimes_{\sZ}
[W_{\g} \mid W_{\h}]_2
\lto
[W_{\f} \mid W_{\h}]_2
$$
to the composition of $\Hom$
for $\Z$-torsors
$\Dim ({\ad_x}(\f))$, $\Dim ({\ad_x}(\g))$, $\Dim ({\ad_x}(\h))$
\end{prop}
\proof
From the definition of $[W_{\f} \mid W_{\g}]_2$
it is clear that it is enough to prove the proposition when
$$
\f \subset \g \subset \h  \mbox{.}
$$

We consider an exact sequence of sheaves on a scheme
 $\Spec \hoo_{x,X}$
\begin{equation} \label{pos}
0 \lto  \f   \lto \g  \lto \g / \f  \lto 0 \mbox{.}
\end{equation}

We restrict  exact sequence~(\ref{pos})
to the subscheme $U = \Spec \hoo_{x,X} \setminus x $,
and we apply  the functor of adeles
on the scheme $U$ to this exact sequence.
Since  adelic rings are exact functors
 (see~\cite{B}, \cite{H}),
we have the following exact sequence of complexes
\begin{equation} \label{pi}
0 \lto  {\ad_x}(\f) \lto {\ad_x}(\g) \lto {\ad_x}(\g / \f) \lto 0
\mbox{.}
\end{equation}

The scheme $U$
is one-dimensional, and a sheaf
 $(\g / \f) \mid_U $
has a support on a zero-dimesional subscheme of $U$.
Therefore an adelic complex
${\ad_x}(\g / \f)$
has only the zero component,
which is equal to $H^0 (U, \g / \f \mid_U )$.
Since
$$
H^0 (U, \g / \f \mid_U )
 =
\mathop{\mathop{\lim}\limits_{\lto}}\limits_n
\Hom\nolimits_{\hat{\oo}_{x,X}}
(m_x^n \; , \; \g / \f ) \mbox{,} $$
we have
\begin{equation}  \label{yakuti}
 \Dim({{\ad_x}(\g / \f)}) = \Dim (H^0 (U, (\g / \f) \mid_U )) =
[\f \mid \g]_2  \mbox{.}
\end{equation}

We consider the long exact cohomological sequence
which is associated with sequence of complexes~(\ref{pi})
$$
0 \to H^0(U, \f \mid_U) \to   H^0(U, \g \mid_U)
\to  H^0 (U, (\g / \f) \mid_U )
\to  H^1(U, \f \mid_U)
\to H^1(U, \g \mid_U)
\to 0  \mbox{.}
$$

We split successivly this sequence on short exact sequences
\begin{equation} \label{p1}
0 \lto \f_{(0)} \lto \g_{(0)} \lto (\g / \f)_{(0)} \lto 0
\end{equation}
\begin{equation}  \label{p2}
0 \lto (\g / \f)_{(0)}  \lto  (\g / \f)_{(1)} \lto  (\g / \f)_{(2)}
\lto 0
\end{equation}
\begin{equation}  \label{p3}
0 \lto   (\g / \f)_{(2)}  \lto \f_{(2)} \lto \g_{(2)} \lto 0  \mbox{.}
\end{equation}
Here $\f_{(0)} = H^0(U, \f \mid_U)$,
$\g_{(0)} = H^0(U, \g \mid_U)$,
$(\g / \f)_{(1)} = H^0 (U, (\g / \f) \mid_U ) $, \\
$\f_{(2)} = H^1(U, \f \mid_U)$,
$\g_{(2)} = H^1(U, \g \mid_U)$.

From the explicit proof of proposition~(\ref{proposi}) it follows
that exact sequences~(\ref{p1})-(\ref{p3})
are exact triples of linearly locally compact
$k$-spaces
with closed kernels, with the induced topology on the kernels,
and with the factor topology on the factorspaces.
Therefore from proposition~\ref{Kap}
we obtain an isomorphism of
 $\Z$-torsors
$$
\Dim(\f_{(0)}) \otimes_{\sZ}
\Dim((\g / \f)_{(0)}) \lto
\Dim(\g_{(0)})
$$
$$
\Dim((\g / \f)_{(0)})
\otimes_{\sZ}
\Dim((\g / \f)_{(2)})
\lto  \Dim((\g / \f)_{(1)})
$$
$$
\Dim((\g / \f)_{(2)}) \otimes_{\sZ}
\Dim(\g_{(2)})
 \lto \Dim(\f_{(2)}) \mbox{.}
$$
Hence and in view of~(\ref{yakuti}) we get
isomorphism~(\ref{ttt}).

The compatibility of isomorphism~(\ref{ttt}) with the filtration
$\f \subset \g \subset \h$
follows from the exactness of the following diagram
of linearly locally compact $k$-spaces:
$$
\begin{array}{ccccccccc}
& & 0 && 0 && 0 && \\
&& \downarrow && \downarrow && \downarrow && \\
0  & \lto & (\g /\f)_{(0)} &   \lto &  (\h / \f)_{(0)} &  \lto   &
(\h / \g)_{(0)} &  \lto  &  0 \\
&& \downarrow && \downarrow && \downarrow && \\
0  &  \lto &  (\g /\f)_{(1)} &   \lto &  (\h / \f)_{(1)} &  \lto &
(\h / \g)_{(1)} &  \lto &  0 \\
&& \downarrow && \downarrow && \downarrow && \\
0 & \lto &  (\g /\f)_{(2)}  &  \lto &  (\h / \f)_{(2)} &  \lto &
(\h / \g)_{(2)} &  \lto &  0 \\
&& \downarrow && \downarrow && \downarrow && \\
& & 0 && 0 && 0 &&
\end{array}
$$
\begin{flushright}
$  \fbox{}  $
\end{flushright}

 We consider a $k$-space
$$
W = W_{\hoo_{x,X}}   \subset \prod_{F \subset \hoo_{x,X}} K_{x,F}
$$
We consider a group
$$
 H = \Frac ( \hoo_{x,X})^*   \mbox{,}
$$
which acts on $\prod\limits_{F \subset \hoo_{x,X}} K_{x,F}$
in the diagonal way.
Similar to section~\ref{s3.6}
we obtain a central extension
\begin{equation}  \label{ccc}
0 \lto \Z \lto \widehat{H}_W \lto
 H
 \lto 1  \mbox{.}
\end{equation}

\begin{prop}  \label{ras}
The central extension~(\ref{ccc}) splits.
\end{prop}
\proof
There is a natural action of the group $H$
on the set of complexes
 $\ad_x (\f)$:
$$
h \in H  \quad : \quad \ad_{x}(\f)  \lto \ad_{x}(h \f)
$$
It follows from the proof of proposition~\ref{proposi}
that this action induces a well-defined action on
 $\Z$-torsors:
\begin{equation} \label{dej}
h \in H  \quad : \quad
\Dim(\ad_{x}(\f))  \lto  \Dim(\ad_{x}(h \f)) \mbox{.}
\end{equation}
This action will be compatible with isomorphism~(\ref{ttt}),
i.e., the following diagram is commutative
$$
\begin{array}{ccc}
[W_{\f} \mid W_{\g}] & \lto &
\Hom\nolimits_{\sdz} ( \Dim(\ad_x (\f)) , \Dim(\ad_x (\g))) \\
\downarrow & & \downarrow \\
{[h\f \mid h\g]}_{2}

&
\lto
& \Hom\nolimits_{\sdz} ( \Dim(\ad_x (h\f)) , \Dim(\ad_x (h\g)))
 \mbox{.} \\
\end{array}
$$
To prove the splitting of central extension~(\ref{ccc})
we define another central extension
$\widehat{H}'_W $ of the group $H$ by  $\Z$
as the set of pairs $(h, d)$,
where $h \in H$ and the element $d$
is from a
 $\Z$-torsor
$\Hom\nolimits_{\sdz} ( \Dim(\ad_x (\hoo_{x,X})) ,
\Dim(\ad_x (h \hoo_{x,X})))$.
The multiplication in
$\widehat{H}'_W $
given by a rule
$(h,d_1) (g, d_2) = (hg, d_1 \otimes h (d_2))$
makes a group from $\widehat{H}'_W$.

Now an isomorphism~(\ref{ttt}) of  $\Z$-torsors
$$
[W \mid h W]_2 \lto
\Hom\nolimits_{\sdz} ( \Dim(\ad_x (\hoo_{x,X})) ,
\Dim(\ad_x (h \hoo_{x,X})))
$$
induces an isomorphism of central extensions:
$\widehat{H}_W \lto \widehat{H}'_W $.
But the central extension $\widehat{H}'_W$
has a canonical splitting,
which maps an element $h \in H$
to the action~(\ref{dej})
of element $h$ on $\Dim (\ad_x (\hoo_{x,X}))$,
i.e., to an element from
$
\Hom\nolimits_{\sdz} ( \Dim(\ad_x (\hoo_{x,X})) ,
\Dim(\ad_x (h \hoo_{x,X})))
$.
\begin{flushright}
$  \fbox{}  $
\end{flushright}

\begin{Th}  \label{te3}
Let $f,g \in \Frac{\hoo_{x,X}}^*$.
Then the following sum contains only a finite number of
non-zero terms and
\begin{equation} \label{summa}
\sum_{F \subset \hoo_{x,X}} \nu_{{x,F}} (f,g) = 0   \mbox{,}
\end{equation}
where the sum is taken over all prime ideals of codimension $1$
in the ring $\hoo_{x,X}$.
\end{Th}
\proof
According to theorem~\ref{te1} we have
$$
\nu_{{x,F}} (f,g) = <f,g>_{K_{x,F}} \mbox{.}
$$
Therefore we will prove~(\ref{summa}) for $<\;,\;>_{K_{x,F}}$.
We choose free
 $\hoo_{x,X}$-modules $\h_1 \subset \hoo_{x,X}$
and $\h_2 \subset  \hoo_{x,X} $
such that
$\h_1  \subset \hoo_{x,X} $,
$\h_2  \supset \hoo_{x,X} $,
$\h_1 \subset f \hoo_{x,X} \subset  \h_2 $,
$\h_1 \subset  g \hoo_{x,X}  \subset  \h_2 $.

Let $I = {\rm Supp \:} (\h_2 / \h_1) $
be a set of prime ideals $F$
which are the support of  $\hoo_{x,X}$-module $\h_2 / \h_1$.
The set $I$ is finite.
(We chose the set $I$ such that it contains all the zeroes
and poles of functions $f$ and $g$
in $\hoo_{x,X}$.)
Let
$$ M_1 = \prod_{F \in I} K_{x,F}  \quad \mbox{,} \quad
M_2 = \prod_{F \ne I} K_{x,F}  \quad  \mbox{,} \quad
M = M_1 \times M_2  \mbox{.}
$$
Then according to proposition~\ref{predl7} we have
$$
<f, g>_{M} = <f, g>_{M_1}  + <f, g>_{M_2}  \mbox{.}
$$

We have
$$
f \prod_{F \ne I} \oo_{K_{x,F}} = \prod_{F \ne I} \oo_{K_{x,F}}
\qquad
\mbox{and}
\qquad
g \prod_{F \ne I} \oo_{K_{x,F}} = \prod_{F \ne I} \oo_{K_{x,F}}
\mbox{,}
$$
Therefore
similar to the proof of item~\ref{pun} of proposition~\ref{predlo}
we get that the central extension $\widehat{H}_{M_2 \cap W}$
splits over a subgroup generated by the elements $f$ and $g$.
Now in consideration of proposition~\ref{kocy}  we get
$$
<f, g>_{M_2} = 0   \mbox{.}
$$

From proposition~\ref{ras}  we get $<f, g>_{M} = 0$.
Therefore we have
\begin{equation} \label{z}
0 = <f,g>_{M_1} = \sum_{F \in I} <f,g>_{K_{x,F}}   \mbox{.}
\end{equation}
For any $F \ne I$ the following formula holds
$$ f \oo_{K_{x,F}} = \oo_{K_{x,F}}  \qquad \mbox{and} \qquad
 g \oo_{K_{x,F}} = \oo_{K_{x,F}}  \mbox{,}
$$
therefore $<f,g>_{K_{x,F}} = 0 $ when $F \ne I$.
The last expression together with equality~(\ref{z})
is equivalent to the statement of theorem.
\begin{flushright}
$  \fbox{}  $
\end{flushright}

\begin{cons}
 Let $f,g \in k(X)^*$.
Then the following sum contains only  a finite number of
non-zero terms and
\begin{equation} \label{summa}
\sum_{C \ni x} \nu_{{x,C}} (f,g) = 0   \mbox{,}
\end{equation}
where the sum is taken over all irreducible curves $C \subset X$
going through the point $x$.
\end{cons}
\proof
Let $J_C \subset \oo_{x,X}$
be an ideal of the curve
 $C$ in the local ring  $\oo_{x,X}$ of the point $x$ on $X$.
 Suppose that with  an ideal
 $J_C \cdot \hoo_{x,X}$
in the completed ring
 $\hoo_{x,X}$
are associated  prime ideals
$F_1, \ldots, F_m $,
then
$$ K_{x,C} = \bigoplus\limits_{i=1}^m K_{x,F_i} \qquad  \mbox{,} \qquad
\nu_{{x,C}} = \bigoplus\limits_{i=1}^m  \nu_{{x,F_i}} \mbox{.}
$$

Suppose that a prime ideal
$F$
of codimension
 $1$ in the ring $\hoo_{x,X}$
are not associated with some ideal
 $J_C \cdot \hoo_{x,X}$
for some irreducible curve
 $C \subset X$,
then
$F = F' \cap \hoo_{x,X}$
for some prime ideal
 $F'$
from a one-dimensional ring $k(X) \cdot \hoo_{x,X}$.
Then we have $<f,g>_{K_{x,F}} = 0$ for $f, g \in k(X)^*$.

Now we apply theorem~\ref{te3}.
\begin{flushright}
$  \fbox{}  $
\end{flushright}

\begin{nt} {\em
The results of this section hold also for singular points $x$
on the  algebraic surface $X$.
As before, we have to consider adelic complexes on the schemes
$U = \Spec \hoo_{x,X} \setminus x$.
Then the key tool is proposition~\ref{proposi}.
For its proof  we have to apply the following lemma
\begin{lemma}   \label{lk}
Let $x \in X$
be any point  on algebraic surface. (The point $x$ can be singular.)
Then for any finitely generated
$\hat{\oo}_{x,X}$-module
$\g$  $k$-spaces
$$
{\rm Ext}^i_{\hat{\oo}_{x,X}}
(\hat{\oo}_{x,X} / m_x^n \, , \, \g)
$$
are finite dimensional.
\end{lemma}
\proof
At first, we compute
${\rm Ext}^i_{\hat{\oo}_{x,X}} (\hat{\oo}_{x,X} / m_x^n \, , \cdot)$
as a functor of the second argument through
an injective resolution of  $\hat{\oo}_{x,X}$-modules.
Hence we get that the  groups   \\
$
{\rm Ext}^i_{\hat{\oo}_{x,X}} (\hat{\oo}_{x,X} / m_x^n \, , \cdot)$
are annuled
by $m_x^n$ as $\oo_{x,X}$-modules.
Now
we compute the groups
${\rm Ext}^i_{\hat{\oo}_{x,X}} (\cdot , \, \g)$
as a functor of the first  argument
through  a free resolution of $\hat{\oo}_{x,X}$-modules.
Hence we get that the groups
${\rm Ext}^i_{\hat{\oo}_{x,X}} (\cdot , \, \g)$
are finitely generated as $\oo_{x,X}$-modules.
Therefore the  $k$-spaces
$
{\rm Ext}^i_{\hat{\oo}_{x,X}}
(\hat{\oo}_{x,X} / m_x^n \, , \, \g)
$
are finite dimensional over $k$.
\begin{flushright}
$  \fbox{}  $
\end{flushright}
}
\end{nt}

\section{Reciprosity laws along the curves.}
\subsection{Adelic complex connected with a curve on a surface}
Let $X$
be a smooth algebraic surface over a field  $k$.
Let $C$
be an irreducible  curve on $X$.

Let $\f$
be a coherent sheaf on $X$.
We construct an adelic complex
 $\da_C(\f)$
of abelian groups in the following way:
$$
\ad_C(\f) =
\mathop{\Lim_{\rightarrow}}_n
\mathop{\Lim_{\leftarrow}}_{m \ge n}
\ad_X (\f \otimes_{\oo_X} J_C^n / J_C^m)  \mbox{.}
$$
Here $J_C$
is an ideal sheaf of the curve $C$  on $X$,
and $\ad_X$
is a functor from coherent sheaves on $X$ to adelic complexes on $X$
 (see~\cite{B}, \cite{H}).

The complex $\da_C(\f)$
has not more than two members,
since the sheaf
 $\f \otimes_{\oo_X} J_C^n / J_C^m$
comes from some infinitesimal neighbourhoods of the curve $C$ in $X$.

\begin{prop}  \label{po}
Let
$$
0 \lto \f \lto \g \lto \h \lto 0
$$
be an exact triple of coherent sheaves on $X$.
Then the following triple of complexes of abelian groups is exact
$$
0 \lto  \ad_C(\f) \lto \ad_C(\g)  \lto \ad_C(\h) \lto 0 \mbox{.}
$$
\end{prop}
\proof
Injective limits preserve the exactness of complexes.
Therefore for any integer $n$
it is enough to prove the exactness of a functor
$$
\mathop{\Lim_{\leftarrow}}_{m \ge n }
\ad_X (\f \otimes_{\oo_X} J_C^n / J_C^m)  \mbox{.}
$$
For any coherent sheaf $\f$ on  $X$ and integer $n$
we define $\f_n = \f \otimes_{\oo_X} J_C^n $.

For any natural $l$
is exact the following sequence of sheaves on $X$
$$
0 \lto   \f_n / J_C^l \g_n \cap \f_n
\lto \g_n /J_C^l \g_n   \lto \h_n/ J_C^l \h_n \lto 0  \mbox{.}
$$

The functor $\ad_X$
is an exact functor, therefore the following sequence of
complexes is exact:
$$
0 \lto  \ad_X (\f_n / J_C^l \g_n \cap \f_n)
\lto
\ad_X(\g_n /J_C^l \g_n)   \lto
\ad_X(\h_n/ J_C^l \h_n) \lto 0  \mbox{.}
$$
For any $l \ge l'$  the maps
$$
\f_n / J_C^l \g_n \cap \f_n
\lto
\f_n / J_C^{l'}  \g_n \cap \f_n
$$
are surjective,
therefore are surjective the following  maps
$$
\ad_X(\f_n / J_C^l \g_n \cap \f_n)
\lto
\ad_X( \f_n / J_C^{l'}  \g_n \cap \f_n ) \mbox{.}
$$
Therefore  the functor  $\mathop{\Lim}\limits_{\leftarrow}$ is exact.
Therefore the following sequence is exact:
\begin{equation}  \label{main}
0 \lto
\mathop{\Lim_{\leftarrow}}_{l \ge 0 }
\ad_X (\f_n / J_C^l \g_n \cap \f_n)  \lto
\mathop{\Lim_{\leftarrow}}_{l \ge 0 }
\ad_X(\g_n /J_C^l \g_n)   \lto
\mathop{\Lim_{\leftarrow}}_{l \ge 0 }
\ad_X(\h_n/ J_C^l \h_n) \lto 0  \mbox{.}
\end{equation}

From the Artin-Rees lemma (see~\cite{AM}) it follows that
there exist a natural number $k$ such that
$
J_C^l \g_n \cap \f_n = J_C^{l-k} (J_C^k \g_n  \cap \f_n)
$
for all $l \ge k$.
Therefore the following maps are well defined:
$$
\f_n / J_C^l \g_n \cap \f_n
\lto
\f_n / J_C^{l-k} \f_n
\qquad \mbox{and} \qquad
\f_n / J_C^l \f_n  \lto \f_n / J_C^l \g_n \cap \f_n  \mbox{.}
$$
Therefore $\ad_X (\f_n / J_C^l \g_n \cap \f_n)$ and
$\ad_X (\f_n / J_C^l \f_n)$
are cofinale projective systems
when $l$ run over natural numbers.
Hence
$$
\mathop{\Lim_{\leftarrow}}_{m \ge n }
\ad_X (\f \otimes_{\oo_X} J_C^n / J_C^m)  =
\mathop{\Lim_{\leftarrow}}_{l \ge 0 }
\ad_X (\f_n / J_C^l \f_n ) =
\mathop{\Lim_{\leftarrow}}_{l \ge 0 }
\ad_X (\f_n / J_C^l \g_n \cap \f_n)
 \mbox{.}
$$
Hence, from~(\ref{main})
and after the passing to injective limit
we get that the functor $\ad_C$ is exact.
\begin{flushright}
$  \fbox{}  $
\end{flushright}

For any coherent sheaf $\f \subset $ on $X$ we have that
the complex $\ad_C (\f)$
has the following form:
$$
 A(\f) \times \prod_{x \in C} B_{x}(\f) \lto D(\f)
$$
From proposition~\ref{po} and similar to the  proof of
proposition~(10.13) from~\cite{AM} we get that
$B_{x} (\f) = B_x (\oo_{X}) \otimes_{\hoo_{x,X}} \f$.
We denote $B_x = B_x (\oo_{X})$.

Let $\f \subset k(X)$
be an invetible sheaf on $X$.
Then the complex $\ad_C(\f)$
has the following explicit form
$$
 \widehat{k(X)}_C
\;
\times
\;
\prod_{x \in C}
(B_x \otimes_{\hoo_{x,X}} \f)
\lto \da_C \mbox{,}
$$
where the product is taken over all points of the curve
$C$,
the field $\widehat{k(X)}_'$
is a completion of the field $k(X)$
with respect to the discrete valuation
given by the irreducible curve $C$,
the ring $B_x$
is an  $\hoo_{x,X}$-submodule in the ring  $K_{x,C}$.

The ring $\da_C$
has the following explicit form
$$
\da_C =  \{f_x \} \in \prod_{x \in C} K_{x,C} \qquad  \mbox{such that}
$$
 for some local parameter $u_C$
which gives the curve $C$
on open $ U \subset X$, for almost all $x \in  C \cap U$
from decomposition
$f_{x} = \sum\limits_i a_{x,i} u_C^i$
we have that for every
 $i$ a collection $\{a_{x,i} \in \widehat{k(C)}_x \}$
is a usual adel on the curve $C$.
(It means that for every fixed $i$
for almost all points $x \in C$ we have
 $a_{x,i} \in \oo_{\widehat{k(C)}_x}$.)

Suppose that the point $x$ is a smooth point on the curve $C$,
then after the choice of local parameters
we have $K_{x,C} = k'((t))((u_C))$ and $B_{x} = k'[[t]]((u_'))$.
($k' = k(x)$ is the residue field of the point $x$.)

For any invertible sheaf $\f \subset k(X)$
we denote
$$
W_{\f} =
\prod_{x \in C}  B_x \otimes_{\hoo_{x,X}} \f  \quad \subset \quad \da_C
\mbox{.}
$$

For any invertible sheaves $\f$ and $\h$
such that $\h \subset \f \subset k(X)$ we have that a
$k$-space
$$ W_{\f} / W_{\h} = \prod_{x \in C}  B_x \otimes_{\hoo_{x,X}}
(\f / \h) $$
is a linearly locally compact space.
Therefore for any invertible sheaves   $\f$, $\g$ from $k(X)$
 it makes  sence  the following definition from section~\ref{razdel}.
$$
[W_{\f} \mid W_{\g}]_2 =
\quad \mathop{\Lim_{\longleftarrow}}_{W_{\h}} \quad  [W_{\f} \mid W_{\g} \;
; \; W_{\h}]_2
\mbox{,}
$$
where the limit is taken over all invertible sheaves $\h \subset \g$,
$\h \subset \f$.

We consider a central extension
\begin{equation} \label{centr}
0 \lto \Z \lto \widehat{k(X)^*}_{W_{\oo_X}} \lto k(X)^*  \lto 1 \mbox{,}
\end{equation}
where $\widehat{k(X)^*}_{W_{\oo_X}}$
consists of a set of pairs
$(g,d)$, $g \in k(X)^*$, $d \in [W_{\oo_X}, g W_{\oo_X}]_2$.
We define a standart  multiplication
$(g_1, d_1) (g_2, d_2) = (g_1 g_2, d_1 g_1(d_2))$.

\begin{prop} \label{ra}
Let the curve $C$ be projective.
Then the central extension~(\ref{centr}) splits.
\end{prop}
\proof
For the invertible sheaf $\f \subset k(X)$
we consider the cohomology of the complex $\ad_C(\f)$.
The cohomology groups $H^i(X, \f \otimes J_C^n / J_C^m)$
are finite dimensional and coincide with the cohomology of the complex
$\ad_X(\f \otimes J_C^n / J_C^m)$.
Therefore the cohomology of the complex $\ad_C(\f)$
coincides with
$$
\mathop{\Lim_{\rightarrow}}_n
\mathop{\Lim_{\leftarrow}}_{m \ge n}
H^i (X, \f \otimes J_C^n / J_C^m) \mbox.
$$
Hence we have that $H^0(\ad(\f))$
is a linearly locally compact $k$-space
with a topology given by inductive and projective limits.

 For $k \le n \le m$  let
$\phi_{knm}: H^1 (X, \f \otimes J_C^n / J_C^m)  \to
H^1 (X, \f \otimes J_C^k / J_C^m) $    be
a natural map.
Let $\Ker \phi_{nm} =
\mathop{\Lim\limits_{\rightarrow}}\limits_k \Ker \phi_{knm} $.
We denote
$$
{H^1}'(\ad_C(\f))=
\mathop{\Lim_{\rightarrow}}_n
\mathop{\Lim_{\leftarrow}}_{m \ge n}
H^i (X, \f \otimes J_C^n / J_C^m) / \Ker \phi_{nm} \mbox.
$$
Then a $k$-space ${H^1}'(\ad_C(\f))$
is a linearly locally compact space,
and we have a natural map
$ {H^1}(\ad_C(\f)) \to  {H^1}'(\ad_C(\f))$.

Here we go from the non-Hausdorff space to the Hausdorff space,
i.e., we take a factorspace by the closure of zero in ${H^1}(\ad_C(\f))$.

Let $\f \subset \g \subset k(X)$ be invertible sheaves,
then we have an exact sequence
\begin{equation} \label{nonstandart}
0 \to H^0 (\ad_C(\f)) \to H^0 (\ad_C(\g)) \to H^0 (\ad_C (\g / \f))
\to {H^1}' (\ad_C(\f)) \to {H^1}' (\ad_C(\g)) \to 0 \mbox{,}
\end{equation}
which one obtains from the usual long exact cohomological
sequence by means of factorspace by the closure of zero in each term.

We denote $\Dim{\ad_C(\f)} = \Hom_{\sZ} (\Dim({H^1}' (\ad_C(\f)) ),
\Dim({H^0} (\ad_C(\f))) )$.
Now by the similar reasons as in the proof of proposition~\ref{stand},
and from exact sequence~(\ref{nonstandart})
we have a canonical isomorphism
\begin{equation} \label{uch}
[W_{\f} \mid W_{\g}]_2
\lto
\Hom\nolimits_{\sZ} ( \Dim ({\ad_C}(\f)), \Dim ({\ad_C}(\g)) )
\end{equation}
for any invertible sheaves
 $\f$ and $\g$ from $k(X)$.

Any element  $h \in k(X)^*$  maps  $\Dim ({\ad_C}(\f))$
to $\Dim ({\ad_C}(h\f))$.
Therefore by the similar reasons as in the proof of proposition~(\ref{ras}),
and from isomorphism~(\ref{uch})
we get that the central extension~(\ref{centr}) splits.
\begin{flushright}
$  \fbox{}  $
\end{flushright}

\subsection{Connection of central extensions,
which come  from a curve and from
a point.}

Let $X$
be a smooth algebraic surface,
$C$
be an irreducible curve on $X$,
$x$ be a point on $C$.
We consider a group $H = \Frac(\hoo_{x,X})^*$.

We consider a central extension
\begin{equation}  \label{first}
0 \lto \Z \lto \widehat{H}_{B_x} \lto H \lto 0 \mbox{,}
\end{equation}
where $\widehat{H}_{B_x}$
is a set of pairs $(g,d)$,
$g \in H$, $d \in [B_x \mid gB_x ]$
with a standart multiplication
$(g_1, d_1)(g_2, d_2) = (g_1 g_2, d_1 g_1(d_2))$.

We consider another central extension
\begin{equation} \label{second}
0 \lto \Z \lto \widehat{H}_{\oo_{K_{x,C}}} \lto H \lto 0 \mbox{,}
\end{equation}
where $\widehat{H}_{\oo_{K_{x,C}}}$
is a set of pairs $(g,d)$,
$g \in H$, $d \in [\oo_{K_{x,C}} \mid g \oo_{K_{x,C}}]$
with a standart multiplication
$(g_1, d_1)(g_2, d_2) = (g_1 g_2, d_1 g_1(d_2))$.
\begin{prop} \label{end}
The central extension~(\ref{first})
is dual to the central extension~(\ref{second}).
\end{prop}
\proof
Let $\f$ be an invertible sheaf, $\f \subset k(X)$.
A scheme $\Spec \hoo_{x,X} \setminus C$ is affine,
therefore by theorem~2 from~\cite{Os}
the cohomology of complex $\ad_{x,C} (\f)$
$$
( B_x \otimes_{\hoo_{x,X}} \f) \,  \times \, (\oo_{K_{x,C}}
\otimes_{\hoo_{x,X}} \f)
\lto K_{x,C}
$$
coincide with cohomology groups $H^i (U, \f \mid U)$,
where $U = \Spec \hoo_{x,X} \setminus x$ is a one-dimensional scheme.
Indeed, the complex $ \ad_{x,C} (\f) $
is a complex from the Krichever correspondence
(or restricted adelic complex) for the one-dimensional scheme $U$,
the closed point $C \cap U$ on the scheme $U$
and the affine scheme $U \setminus (C \cap U)$.

From the proof of proposition~\ref{proposi} it follows
that a $\Z$-torsor
$\Dim(\ad_{x,C} (\f)) = \Hom_{\sZ}
(\Dim(H^1 (\ad_{x,C} (\f))), \Dim(H^0 (\ad_{x,C} (\f))))$
is well-defined.
Let $\f \subset \g$,
then we apply the functor $\ad_{x,C}$,
and from the long cohomological sequence
we obtain the existence of the following isomorphism
$$
[B_x \otimes_{\hoo_{x,X}} \f \mid  B_x \otimes_{\hoo_{x,X}} \g]_2
\otimes_{\sZ}
  [\oo_{K_{x,C}} \otimes_{\hoo_{x,X}} \f
\mid   \oo_{K_{x,C}} \otimes_{\hoo_{x,X}} \g]_2
\lto   \qquad \qquad \qquad \qquad \qquad
$$
\begin{equation} \label{last}
\qquad \qquad \qquad \qquad \qquad \qquad
\qquad   \lto
\Hom\nolimits_{\sZ} ( \Dim (\ad_{x,C}(\f)), \Dim (\ad_{x,C}(\g)))
\end{equation}
We consider a central extension $\hat{H}$
which consists of a set of pairs $(h, d)$,
where $h \in H$, $d \in
[B_x  \mid h B_x ]_2
\otimes_{\sZ}
  [\oo_{K_{x,C}}
\mid  h \oo_{K_{x,C}} ]_2
$
with a standard multiplication
$(g_1, d_1)(g_2, d_2) = (g_1g_2, d_1 g_1(d_2))$.
By the similar reasons as in the proof of proposition~\ref{ras}
and from~(\ref{last}) we get that the central extension
 $\hat{H}$
splits over the group $H$.
\begin{flushright}
$  \fbox{}  $
\end{flushright}

\begin{nt} \label{zamechan} {\em
From this proposition and from~(\ref{koc})
it follows that the commutator
of liftings of elements in central extension~(\ref{first})
is equal to minus commutator of
liftings of elements in central extension~(\ref{second}). }
\end{nt}

\begin{Th}
Let $f,g \in k(X)^*$.
Then the following sum contains only a finite number
of non-zero terms and
\begin{equation} \label{summa}
\sum_{x \in  C} \nu_{{x,C}} (f,g) = 0   \mbox{,}
\end{equation}
where the sum is taken over all points $x$
of the  irreducible projective curve
 $'$ on the surface $X$.
\end{Th}
\proof
From proposition~\ref{ra}
and by the similar reasons as in the  proof of theorem~\ref{te3}
we get that the sum of commutators of liftings of elements from $k(X)^*$,
which are computed from central extensions (\ref{first}),
is equal to zero.

Now from remark~\ref{zamechan}  and theorem~\ref{te1}
we get that every commutator of lifting of elements,
which is computed from central extension (\ref{first}),
coincides with $-\nu_{{x,C}}$.
\begin{flushright}
$  \fbox{}  $
\end{flushright}

Steklov Mathematical Institute

d\_osipov@mi.ras.ru

\end{document}